\documentstyle[12pt]{amsart}
\let\text=\mbox
\newcommand{\skipthistext}[1]{}
\newcommand{\C}{{\mathbb C}}
\newcommand{\R}{{\mathbb R}}
\newcommand{\Z}{{\mathbb Z}}
\newcommand{\Q}{{\mathbb Q}}

\newcommand{\I}{{\mathbb I}}

\renewcommand{\pf}{\noindent {\bf Proof} \hspace{2mm}}
              \newcommand{\J}{{\cal J}}
              \newcommand{\M}{{\cal M}}
              
              \newcommand{\U}{{\cal U}}

              \newcommand{\E}{{\cal E}}

              \renewcommand{\O}{{\cal O}}

\newcommand{\abs}[1]{\centerline{\sc Abstract}
\vspace{3mm}
\centerline{\parbox{150mm}{\small #1}}}
\newcounter{thm}

\def\thmnumber{
\addtocounter{thm}{1}\if\arabic{section}0\relax
\else\arabic{section}.\fi%
\if\arabic{subsection}0\relax
\else\arabic{subsection}.\fi%
\arabic{thm}
\hspace{2mm}}

\def\setthmnumber{\setcounter{thm}{0}}

\newenvironment{lem}{
\vspace{3mm}\par\noindent
{\bf Lemma \thmnumber} \it}
{\vspace{3mm}\par\rm}

\newcommand{\sectioni}[1]{\setthmnumber\setcounter{subsection}{0}
\section{#1}}

\renewcommand{\subsection}[1]{%
\vspace{7mm}\par\noindent%
\addtocounter{subsection}{1}%
{\sc\arabic{section}.\arabic{subsection} \ {#1}}
\vspace{4mm}\par\noindent\setthmnumber}

\begin{document}
\title{On a theorem of Peters on automorphisms of K\"{a}hler surfaces}
\author{Weimin Chen}
\date{\today}
\maketitle
\abs{For any K\"{a}hler surface which admits no nonzero holomorphic
vectorfields, we consider the group of holomorphic automorphisms
which induce identity on the second rational cohomology. Assuming the
canonical linear system is without base points and fixed
components, C.A.M. Peters \cite{Pe} showed that this group is
trivial except when the K\"{a}hler surface is of general type and
either $c_1^2=2c_2$ or $c_1^2=3c_2$ holds. Moreover, this group is a
$2$-group in the former case, and is a $3$-group in the latter.
The purpose of this note is to give further information about this
group. In particular, we show that $c_1^2$ is divisible by the
order of the group. Our argument is based on the results of C.H. Taubes
in \cite{T1, T2} on symplectic $4$-manifolds, which are applied here in
an equivariant setting.
}

\sectioni{Introduction}

Let $X$ be a K\"{a}hler surface with $H^0(X,T_X)=0$, ie., $X$
admits no nonzero holomorphic vectorfields, and let
$\text{Aut}(X)^o\subset\text{Aut}(X)$ be the subgroup of
holomorphic automorphisms of $X$ which operates trivially on
$H^2(X;\Q)$. By a result of Lieberman \cite{Lie}, $\text{Aut}(X)^o$
is finite.

In \cite{Pe}, C.A.M. Peters proved the following theorem
concerning $\text{Aut}(X)^o$.

\vspace{2mm}

\noindent{\bf Theorem} (Peters)\hspace{2mm}
{\em Let $X$ be a K\"{a}hler surface with $H^0(X,T_X)=0$. Suppose
$|K_X|$ is without base points and fixed components. Then
$g\in\text{Aut}(X)^o$ is trivial unless $X$ is a surface of
general type and either
\begin{itemize}
\item [{(i)}] $c_1^2=2c_2$, and $|g|$ is a power of $2$, or
\item [{(ii)}] $c_1^2=3c_2$, $|g|$ is a power of $3$ and moreover,
$g$ acts trivially on $H^\ast(X;\Q)$.
\end{itemize}
}

\vspace{2mm}

The purpose of this note is to give further information about
$\text{Aut}(X)^o$ for the two exceptional cases in Peters'
theorem. In particular, we show that $c_1^2$ must be divisible by
the order of $\text{Aut}(X)^o$. Before stating our theorem, we first
have a digression on free actions of a finite group on Riemann surfaces.

Let $G$ be a finite group and $\Sigma_m$ be a Riemann surface of genus $m$
such that $G$ acts freely on $\Sigma_m$ via orientation-preserving
homeomorphisms, and let $\Sigma_n\equiv\Sigma_m/G$ be the quotient Riemann
surface which has genus $n$. Then the following are easily seen:
\begin{itemize}
\item [{(a)}] $\chi(\Sigma_m)=|G|\cdot\chi(\Sigma_n)$, or equivalently,
$m-1=|G|\cdot (n-1)$
\item [{(b)}] $G=\pi_1(\Sigma_n)/\pi_1(\Sigma_m)$ where $\pi_1(\Sigma_m)$
is naturally regarded as a normal subgroup of $\pi_1(\Sigma_n)$ under the
regular covering $\Sigma_m\rightarrow\Sigma_n$
\end{itemize}
With the preceding understood, we introduce the following terminology.
For any finite group $G$, we will call the minimal genus of a Riemann
surface which admits a free $G$-action the {\it free genus} of $G$.
The free genus of a finite group $G$ is closely related to the minimal number
of generators of $G$ (which is simply the rank of $G$ when $G$ is abelian).

\vspace{2mm}

\noindent{\bf Lemma}\hspace{2mm}
{\em
Let $r$ be the minimal number of generators of $G$, and let $[x]$
be the greatest integer less than or equal to $x$. Then
$$
([2^{-1}(r+1)]-1)\cdot |G|+1\leq \text{{\em free genus of }}G
\leq (r-1)\cdot |G|+1.
$$
}
\vspace{2mm}

\pf
The left hand side of the inequality follows easily from the assertions
(a), (b) above. As for the right hand side, we appeal to
the following construction which was pointed out to us by J. McCarthy,
compare also \cite{Fu}.

Let $c_1,\cdots,c_r$ be a set of generators of $G$. Then $G$ can be
realized as a quotient group of
$$
\pi_1(\Sigma_r)=\langle x_1,y_1,\cdots,x_r,y_r\mid
x_1y_1x_1^{-1}y_1^{-1}\cdots x_ry_rx_r^{-1}y_r^{-1}=1\rangle
$$
under the homomorphism $x_i\mapsto c_i$, $y_i\mapsto 1$, $i=1,\cdots,r$.
Let $\Sigma_\delta\rightarrow\Sigma_r$ be the corresponding regular
covering (note that $G$ is finite). Then $(\delta-1)=(r-1)\cdot |G|$,
and the free genus of $G$ is less than or equal to $\delta$. The
lemma follows immediately.

\hfill $\Box$

Now we state our main result.

\vspace{2mm}

\noindent{\bf Theorem}\hspace{2mm}
{\em Let $X$ be a K\"{a}hler surface as in Peters' theorem such that
$\text{Aut}(X)^o$ is nontrivial. Then the following conclusions hold.
\begin{itemize}
\item [{(a)}] Each $g\in\text{Aut}(X)^o$ has order $2$ or $3$.
In particular, in the case of $c_1^2=2c_2$, $\text{Aut}(X)^o$ is an
elementary abelian $2$-group {\em(}ie., product of copies of $\Z_2${\em)}.
\item [{(b)}] $c_1^2$ is divisible by the order of $\text{Aut}(X)^o$.
\item [{(c)}] $c_1^2\geq\max{\em(}\text{{\em free genus of }}
\text{Aut}(X)^o-1, |\text{Aut}(X)^o|{\em)}$.
\end{itemize}
}

\vspace{2mm}

The canonical map of minimal surfaces of general type was
systematically studied by A. Beauville in \cite{Be}, particularly
for the case where the arithmetic genus is relatively large. Based on
Beauville's theorem, J.-X. Cai in \cite{Cai} showed that for a minimal
surface $X$ of general type with $\chi(\O_X)>188$, $\text{Aut}(X)^o$ is
either cyclic of order less than $5$, or is $\Z_2\times\Z_2$. Moreover,
for the case where $|K_X|$ has no base points or fixed components,
it was shown that Beauville's theorem directly implies that when
$\chi(\O_X)\geq 31$, the order of $\text{Aut}(X)^o$ is less than $5$
(cf. \cite{Cai}, page 347). Combined with this observation, we arrived
at the following

\vspace{4mm}

\noindent{\bf Corollary}\hspace{2mm}
{\em Let $X$ be a K\"{a}hler surface as in Peters' theorem with
$\text{Aut}(X)^o$ nontrivial. Then the following must hold.
\begin{itemize}
\item [{(a)}] For the case of $c_1^2=2c_2$, $\text{Aut}(X)^o$ is an
elementary abelian $2$-group of rank $\leq 6$.
\item [{(b)}] For the case of $c_1^2=3c_2$, $\text{Aut}(X)^o$ is a $3$-group
of order $\leq 243$.
\end{itemize}
}

\vspace{2mm}

The proof of our main result, which is given in the next section, is divided
into two parts. In Part 1 we give a proof for part (a) of the theorem, which is
a refined version of Peters' argument in \cite{Pe}, and as in \cite{Pe},
is based on application of $G$-index theorems and the Miyaoka-Yau inequality
$c_1^2\leq 3c_2$. Part 2 is concerned with (b) and (c). The new
ingredient here is the application of the results of C.H. Taubes in
\cite{T1, T2} in an equivariant context. (See also \cite{C, CK}.)
More concretely, we showed that there is a finite set of disjoint, embedded
surfaces $\{C_i\}$ in $X$ such that (1) each $C_i$ lies in the complement
of the exceptional orbits of $\text{Aut}(X)^o$, (2) each $C_i$ is invariant
under the action of $\text{Aut}(X)^o$, and (3) 
$c_1^2=\sum_i (\text{genus}(C_i)-1)$ holds, from which (b) and (c)
of the theorem follow. 

\vspace{3mm}

\centerline{\bf Acknowledgments}

\vspace{3mm}

Benefits from joint work \cite{CK} with S. Kwasik and helpful communications
with J. McCarthy are gratefully acknowledged. The author is partially
supported by NSF grant DMS-0304956.

\sectioni{The proof}

\noindent{\bf Part 1}. We start with the following lemma about the local
representation of $g\in\text{Aut}(X)^o$ at a fixed point, which is Lemma
2 in Peters \cite{Pe}. We wish to point out that the order of $g$ is not
necessarily prime here (which is assumed in \cite{Pe}). We set $\mu_p\equiv
\exp(\frac{2\pi i}{p})$ below.

\begin{lem}
Let $x\in X$ be a fixed point of $1\neq g\in\text{Aut}(X)^o$. Then $x$ is
an isolated fixed point and the action of $g$ near $x$ is given by
$(z_1,z_2)\mapsto (\mu_p^k z_1,\mu_p^{-k}z_2)$ for some integer $k$,
where $p=|g|$, and $k$ is relatively prime to $p$.
\end{lem}

\pf
The proof is the same as in \cite{Pe}. Because $|K_X|$ is without base
points and fixed components, there exists a holomorphic $2$-form $\omega$
such that $\omega(x)\neq 0$. On the other hand, since $g$ acts trivially
on $H^2(X;\Q)$, $g^\ast\omega=\omega$, which implies the lemma.

\hfill $\Box$

The next lemma gives a lower bound on the number of fixed points of a certain
orientation-preserving self-diffeomorphism of a $4$-manifold having only
isolated fixed points (compare Lemma 3 of \cite{Pe}).

\begin{lem}
Let $f:M\rightarrow M$ be an orientation-preserving, periodic,
self-diffeomorphism of a $4$-manifold which has only isolated fixed points.
If $f$ induces identity on $H^2(M;\Q)$, then the number of fixed points of
$f$ is bounded below by the Euler characteristic of $M$.
\end{lem}

\pf
Let $n$ be the number of fixed points of $f$. Then by the Lefschetz
fixed point theorem,
$$
n=\sum_{k=0}^4 (-1)^k \text{Trace}(f^\ast|H^k(M;\Q))=
2+b_2-\text{Trace}(f^\ast|H^1(M;\Q))-\text{Trace}(f^\ast|H^3(M;\Q)),
$$
because $f$ is periodic, orientation-preserving, and induces identity on
$H^2(M;\Q)$. On the other hand, each $H^k(M;\Q)$, $k=1,3$, is decomposed
into $f^\ast$-invariant direct sum $\oplus_i V_i$, where $V_i$ is either
$1$-dimensional with $f^\ast|V_i=\pm 1$, or $2$-dimensional with $f^\ast|V_i$
being a rotation in $V_i$. In any event, $\text{Trace}(f^\ast|V_i)\leq \dim V_i$,
and consequently, $\text{Trace}(f^\ast|H^k(M;\Q))\leq b_k$ for
$k=1,3$. The lemma follows immediately.

\hfill $\Box$

Now we are ready to give a proof for part (a) of the theorem. To this
end, we first recall the following version of the $G$-signature theorem
for a cyclic group action of $G$ on a $4$-manifold $M$, which has only
isolated fixed points (cf. \cite{HZ}, Equation (12) on page 177)
$$
|G|\cdot\text{sign}(M/G)=\text{sign}(M)+\sum_m\text{def}_m.
$$
Here in the above formula, $m\in M$ is running over the set of exceptional
orbits, and $\text{def}_m$ stands for the signature defect. For $m$ with
isotropy subgroup of order $p$ and local representation
$(z_1, z_2)\mapsto (\mu_{p}^k z_1,\mu_p^{kq} z_2)$ where $k,q$
are relatively prime to $p$, the signature defect is given by
the following formula (cf. \cite{HZ}, Equation (19) on page 179)
$$
\text{def}_m=I_{p,q}\equiv\sum_{k=1}^p\frac{(1+\mu_p^k)(1+\mu_p^{kq})}
{(1-\mu_p^k)(1-\mu_p^{kq})}.
$$

The following lemma computes $I_{p,q}$ for the case of $q=-1$ (see also
Lemma 2.3 in \cite{CK}).

\begin{lem}
$I_{p,-1}=\frac{1}{3}(p-1)(p-2).$
\end{lem}

\pf
$I_{p,q}$ can be computed in terms of Dedekind sum $s(q,p)$ (cf. \cite{HZ},
page 92),  where
$$
s(q,p)=\sum_{k=1}^p((\frac{k}{p}))\;((\frac{kq}{p}))
$$
with
$$
((x))=\left\{\begin{array}{cc}
x-[x]-\frac{1}{2} & \mbox{if } x\in\R\setminus\Z\\
0 & \mbox{if } x\in\Z.\\
\end{array} \right.
$$
(Here $[x]$ stands for the greatest integer less than or equal to $x$.)

In fact, Equation (24) in \cite{HZ}, page 180, gives
$$
I_{p,q}=-4p\cdot s(q,p),
$$
with the expression $6p\cdot s(q,p)$ given by
(cf. \cite{HZ}, Equations (10) and (9) on page 94)
$$
6p\cdot s(q,p)=(p-1)(2pq-q-\frac{3p}{2})-6f_p(q),
$$
where $f_p(q)=\sum_{k=1}^{p-1}k[\frac{kq}{p}]$. Since
$f_p(-1)=\sum_{k=1}^{p-1}k\cdot (-1)=\frac{1}{2}(1-p)p$,
one obtains
$$
I_{p,-1}=\frac{1}{3}(p-1)(p-2)
$$
as claimed.

\hfill $\Box$

With the preceding preparation, we now give a proof for part (a) of the
theorem. Suppose to the contrary that
$g\in\text{Aut}(X)^o$ has order other than $2$ or $3$. Then by Peters'
theorem, and by passing to a suitable power of $g$, we may assume without
loss of generality that $|g|=p^2$, where $p=2$ or $3$. Now consider the
cyclic action of $G\equiv \langle g\rangle$ on the K\"{a}hler surface
$X$. Since $G$ operates trivially on $H^2(X;\Q)$, the $G$-signature
theorem implies
$$
(p^2-1)\cdot \text{sign}(X)=\sum_{m\in X^g}\text{def}_m+
\sum_{m\in X^{g^p}}\text{def}_m.
$$
By Lemma 2.1 and Lemma 2.3, $\text{def}_m\geq 0$ for any $m\in X^{g^p}$,
and by Lemma 2.2, $|X^g|\geq c_2$. It follows easily that (cf.
Lemma 2.3)
$$
\frac{1}{3}(p^2-1)(c_1^2-2c_2)\geq \frac{1}{3}(p^2-1)(p^2-2)c_2,
$$
which contradicts the Miyaoka-Yau inequality $c_1^2\leq 3c_2$ (cf. \cite{BPV}).
Hence part (a) of the theorem.

\vspace{3mm}

\noindent{\bf Part 2}. We begin by recalling the relevant theorems
of Taubes in \cite{T1, T2}, which are cast here in an equivariant setting.
See also \cite{C, CK}.

Let $(M,\omega)$ be a symplectic $4$-manifold, and $G$ be a finite group
acting on $(M,\omega)$ via symplectomorphisms. Denote by $b_G^{2,+}$ the
dimension of the maximal subspace of $H^2(M;\Q)$ over which the
cup product is positive and the induced action of $G$ is identity.

For any given $G$-equivariant $\omega$-compatible almost complex structure
$J$, consider the associated $G$-equivariant Riemannian metric
$g=\omega(\cdot,J(\cdot))$. There is a canonical
$G$-$Spin^\C$ structure on $M$ such that the associated $U(2)$
$G$-bundles are $S_{+}^0=\I\oplus K^{-1}$ and $S_{-}^0=T^{0,1}M$,
where $\I$ is the trivial $G$-bundle and $K=\det T^{1,0}M$ is the
canonical bundle. Note that in this setup, the associated
Seiberg-Witten equations for a pair $(A,\psi)$
$$
D_A\psi=0 \mbox{ and }
P_{+}F_A=\frac{1}{4}\tau(\psi\otimes\psi^\ast)+\mu
$$
are $G$-equivariant, where $A$ is a $G$-equivariant
$U(1)$-connection on $\det S_{+}^0$ and $\psi$ is a $G$-equivariant
smooth section of $S_{+}^0$, and $\mu$ is some fixed
$G$-equivariant, imaginary valued self-dual $2$-form.
According to Taubes \cite{T1}, there is a canonical (up to gauge
equivalence) connection $A_0$ on $K^{-1}=\det S_{+}^0$, such that
if we set $u_0\equiv (1,0)\in\Gamma(\I\oplus K^{-1})$, then for
any $r>0$, $(A_0,\sqrt{r}u_0)$ (which is clearly $G$-equivariant) satisfies
the Seiberg-Witten equations with
$\mu=-\frac{i}{4}r\cdot\omega+P_{+}F_{A_0}$. Moreover,
$(A_0,\sqrt{r}u_0)$ is the only solution (up to gauge equivalence)
when $r>0$ is sufficiently large, which is also non-degenerate.
Thus assuming $b_G^{2,+}\geq 2$, if we define a $G$-equivariant
Seiberg-Witten invariant as an algebraic count of the solutions to
the $G$-equivariant Seiberg-Witten equations, Taubes' theorem in \cite{T1}
implies that the corresponding invariant equals $\pm 1$
for the canonical $G$-$Spin^\C$ structure $S_{+}^0\oplus S_{-}^0$.

Consequently, as in the non-equivariant case (assuming $b_G^{2,+}\geq
2$), the $G$-equivariant Seiberg-Witten invariant equals $\pm 1$
for the $G$-$Spin^\C$ structure $K\otimes (S_{+}^0\oplus S_{-}^0)$.
This means that for any $r>0$, there is a solution $(A,\psi)$ to
the $G$-equivariant Seiberg-Witten equations with
$\mu=-\frac{i}{4}r\cdot\omega+P_{+}F_{A_0}$. Particularly, $(A,\psi)$
is $G$-equivariant. Now write $\psi=\sqrt{r}(\alpha,\beta)\in
\Gamma(K\oplus\I)$ (note: $K\otimes S_{+}^0=K\oplus\I$). Then according
to Taubes \cite{T2}, the zero
set $\alpha^{-1}(0)$ converges as $r\rightarrow\infty$ to a finite set
of $J$-holomorphic curves $\{C_i\}$, such that for some integers
$n_i>0$, the canonical class $c_1(K)$ is Poincar\'{e} dual to the
fundamental class of $\sum_i n_i C_i$. The crucial observation here
is that $\cup_i C_i$ is $G$-invariant. This is because $\alpha$ is a
$G$-equivariant section of $K$, and as $r\rightarrow\infty$, $\alpha^{-1}(0)$
converges to $\cup_i C_i$ with respect to the natural distance function
on $M$.

With the preceding understood, the following lemma gives certain regularity
about the $J$-holomorphic curves $\{C_i\}$ for a generic $G$-equivariant $J$,
provided that $(M,\omega)$ and the $G$-action satisfy some further
conditions.

\begin{lem}
Suppose $(M,\omega)$ is minimal, and the action of $G$ is
pseudofree such that for any $m\in M$, the representation
of the isotropy subgroup $G_m$ on the tangent space of $m$ is
contained in $SL_2(\C)$ with respect to a $G$-equivariant,
$\omega$-compatible almost complex structure $J_0$. Then for a generic
$G$-equivariant $J$, the $J$-holomorphic curves $\{C_i\}$ are disjoint,
embedded, all contained in the complement of the exceptional orbits,
with $n_i=1$ for any $C_i$ such that $C_i^2>0$.
\end{lem}

\pf
The non-equivariant version of this result was due to Taubes \cite{T2},
whose proof was based on transversality for moduli spaces of $J$-holomorphic
maps and the adjunction formula. The proof of Lemma 2.4 is a
somewhat equivariant version of that in Taubes \cite{T2}.

Let $H\subset G$ be any subgroup (here $H$ is allowed to be trivial), and
let $\Sigma$ be a Riemann surface which admits a holomorphic $H$-action.
We shall consider the transversality problem for the moduli space of
pseudoholomorphic maps $f:\Sigma\rightarrow M$ where $f$ is equivariant,
ie., $f\circ h=h\circ f$ for any $h\in H$.

To this end, we put the problem in the Fredholm framework as follows. Fix a
sufficiently large $r>0$, we consider the space
$$
[\Sigma;M]^H\equiv \{f:\Sigma\rightarrow M\mid f\text{ is of $C^r$ class,
and } f\circ h=h\circ f, \forall h\in H\},
$$
and the space $\J$ of $G$-equivariant, $\omega$-compatible almost complex
structures $J$ of $C^r$ class, which equals $J_0$ in a fixed neighborhood
of the exceptional orbits. The tangent space $T_f$ of $[\Sigma;M]^H$ at
$f$ is the Banach space of equivariant $C^r$-sections of the $H$-bundle
$f^\ast TM$, and the tangent space $T_J$ of $\J$ at $J$ is the Banach
space of equivariant $C^r$-sections $A$ of the $G$-bundle $\text{End }TM$,
which obeys (1) $AJ+JA=0$, (2) $A^t=A$ (here the transpose $A^t$ is taken
with respect to the metric $\omega(\cdot,J(\cdot))$), and (3)
$A$ vanishes in a fixed neighborhood of the exceptional orbits.

We will also need to consider the moduli space $\M_\Sigma^H$ of
$H$-equivariant complex structures $j$ on $\Sigma$, which is generally
a finite dimensional complex orbifold. For technical reason, we will
cover $\M_\Sigma^H$ by countably many open sets of form $U=\hat{U}/G_U$,
where $\hat{U}$ is a complex manifold and $G_U$ is a finite group, and
work instead with each $\hat{U}$. Denote by $\U$ any one of the $\hat{U}$'s.

With the preceding understood, for any $(f,J,j)\in [\Sigma;M]^H\times
\J\times\U$, let $\E_{(f,J,j)}$ be the Banach space of equivariant
$C^{r-1}$-sections $s$ of the $H$-bundle $\text{Hom}(T\Sigma,f^\ast TM)$
which obeys $s\circ j=-J\circ s$. Then there is a Banach bundle
$\E\rightarrow [\Sigma;M]^H\times\J\times\U$ whose fiber at $(f,J,j)$ is
$\E_{(f,J,j)}$. The zero set of the smooth section $\underline{L}:
[\Sigma;M]^H\times\J\times\U\rightarrow \E$, where
$$
\underline{L}(f,J,j)=df+J\circ df\circ j,
$$
consists of triples $(f,J,j)$ such that $f$ is equivariant (ie.,
$f\circ h=h\circ f, \forall h\in H$), and is $J$-holomorphic with
respect to the complex structure $j$ on $\Sigma$. We consider the subspace
$$
\widetilde{\M}_{H,\Sigma,\U}\equiv \{(f,J,j)\in\underline{L}^{-1}(0)\mid
\text{$f$ is nonconstant and not multiply covered}\},
$$
and state the promised transversality result in the following

\vspace{2mm}

\noindent{\bf Claim 1}\hspace{2mm}
{\em The subspace $\widetilde{\M}_{H,\Sigma,\U}\subset
[\Sigma;M]^H\times\J\times\U$ is a Banach submanifold.
}

\vspace{2mm}

The proof of Claim 1 is based on the same analysis as in the non-equivariant
setting (cf. Proposition 3.2.1 in \cite{McDS}). Recall that the key point in
the argument is that one is able to construct an $A\in T_J$ to kill
the cokernel of $D\underline{L}_{(f,J,j)}|T_f$, where
$D\underline{L}_{(f,J,j)}$ is the linearization of $\underline{L}$ at
$(f,J,j)$. In the present, equivariant context, one needs to construct
an $A\in T_J$ where there is an additional requirement that $A$ must be
$G$-equivariant. But it is easily seen that this can be done, because
for an open, dense subset of $\Sigma$, its image under $f$ lies in the
complement of the exceptional orbits, where the action of $G$ is free.

Now consider the projection $\pi:\widetilde{\M}_{H,\Sigma,\U}
\rightarrow\J$, which is a Fredholm map between Banach manifolds.
By the Sard-Smale theorem \cite{Sm}, there is a Baire set $\J_{reg}
\subset\J$, such that for $J\in\J_{reg}$, the differential $d\pi$ is
surjective along $\pi^{-1}(J)$. For any $J\in\J_{reg}$, we set
$$
\widetilde{\M}_{H,\Sigma,\U}^J\equiv\pi^{-1}(J)=\{(f,j)\mid (f,J,j)
\in\underline{L}^{-1}(0)\},
$$
which is a finite dimensional manifold when nonempty.

Since the set of data $\{(H,\Sigma,\U)\}$ is countable, it follows easily
that there is a Baire set $\J_0=\cap\J_{reg}$, such that for any
$J\in\J_0$, which particularly may be chosen smooth,
$\widetilde{\M}_{H,\Sigma,\U}^J$ is a finite dimensional manifold when
nonempty for all $(H,\Sigma,\U)$.

The next crucial step in the argument is to compute the dimension of
$\widetilde{\M}_{H,\Sigma,\U}^J$. To this end, we observe that the
dimension of $\widetilde{\M}_{H,\Sigma,\U}^J$ at $(f,j)$ is given by
the sum of the index of $D\underline{L}_{(f,J,j)}|T_f$ with the (real)
dimension of the moduli space $\M_\Sigma^H$ of $H$-equivariant complex
structures
on $\Sigma$. The index of $D\underline{L}_{(f,J,j)}|T_f$ can be computed
by the Riemann-Roch theorem for orbit spaces in Atiyah-Singer \cite{AS},
or more generally, the index formula for Cauchy-Riemann type operators
over orbifold Riemann surfaces, cf. Lemma 3.2.4 of \cite{CR}.

We begin by introducing some notations. Given any $(f,j)\in
\widetilde{\M}_{H,\Sigma,\U}^J$, we pick for each exceptional orbit in
$\Sigma$ a point $z_i$ from it, where $i=1,2,\cdots,k$.
We denote by $m_i$ the order of isotropy at $z_i$,
and to each $z_i$, we assign a pair of {\it rotation numbers} $(m_{i,1},
m_{i,2})$ with $0<m_{i,1}, m_{i,2}<m_i$ as follows: let $h_i\in H$ be the
unique element in the isotropy subgroup at $z_i$ whose action near $z_i$
is given by a counterclockwise rotation of angle $\frac{2\pi}{m_i}$, then
the action of $h_i$ on the tangent space of $p_i\equiv f(z_i)$ is given by
$(\xi_1,\xi_2)\mapsto (\mu_{m_i}^{m_{i,1}}\xi_1,\mu_{m_i}^{m_{i,2}}\xi_2)$.
Clearly, the rotation numbers $(m_{i,1},m_{i,2})$ depend only on the
exceptional orbit which $z_i$ lies in. Finally, we set
$\Gamma\equiv\Sigma/H$, which is an orbifold Riemann surface with
orbifold points $w_i=[z_i]$ of orders $m_i$, $i=1,2,\cdots,k$. We denote
by $g_{|\Gamma|}$ the genus of the underlying surface of $\Gamma$.

With these notations, the index of $D\underline{L}_{(f,J,j)}|T_f$ is given
by $2d_{(f,j)}$ where $d_{(f,j)}\in\Z$ and
$$
d_{(f,j)}=\frac{1}{|H|}c_1(TM)\cdot f_\ast([\Sigma])+2-2g_{|\Gamma|}-
\sum_{i=1}^k\frac{m_{i,1}+m_{i,2}}{m_i}.
$$
On the other hand, the moduli space $\M_\Sigma^H$ of $H$-equivariant
complex structures on $\Sigma$ can be identified with the moduli space
of complex structures on the marked Riemann surface $(\Gamma,\{w_i\})$.
Thus we have
$$
\dim_\C \M_\Sigma^H=\left\{\begin{array}{ll}
0 & \text{if } g_{|\Gamma|}=0, k\leq 3,\\
k-3 &  \text{if } g_{|\Gamma|}=0, k>3,\\
1 &  \text{if } g_{|\Gamma|}=1, k=0,\\
k-1 & \text{if } g_{|\Gamma|}=1, k>0,\\
3g_{|\Gamma|}-3+k & \text{if } g_{|\Gamma|}\geq 2.\\
\end{array} \right.
$$

Now here is the crucial consequence of the assumption we made in
Lemma 2.4 that for any $m\in M$, the representation of the isotropy
subgroup $G_m$ on the tangent space of $m$ is contained in $SL_2(\C)$
with respect to a $G$-equivariant, $\omega$-compatible almost complex
structure $J_0$: for any $i=1,2,\cdots,k$, the rotation numbers
$(m_{i,1},m_{i,2})$ obey $m_{i,1}+m_{i,2}=m_i$.

Now assuming $\widetilde{\M}_{H,\Sigma,\U}^J\neq\emptyset$, we see
that $d_{(f,j)}+\dim_\C \M_\Sigma^H\geq 0$ must hold for any
$(f,j)\in\widetilde{\M}_{H,\Sigma,\U}^J$. With $m_{i,1}+m_{i,2}=m_i$,
this gives
$$
\begin{array}{ll}
{\em (i)} \hspace{2mm}|H|^{-1}c_1(TM)\cdot f_\ast([\Sigma])+2-k\geq 0 &
\text{if } g_{|\Gamma|}=0, k\leq 3,\\
{\em (ii)} \hspace{2mm}|H|^{-1}c_1(TM)\cdot f_\ast([\Sigma])-1\geq 0 &
\text{if } g_{|\Gamma|}=0, k>3,\\
{\em (iii)}  \hspace{2mm}|H|^{-1}c_1(TM)\cdot f_\ast([\Sigma])+1\geq 0 &
\text{if } g_{|\Gamma|}=1, k=0,\\
{\em (iv)}  \hspace{2mm}|H|^{-1}c_1(TM)\cdot f_\ast([\Sigma])-1\geq 0 &
\text{if } g_{|\Gamma|}=1, k>0,\\
{\em (v)}  \hspace{2mm}|H|^{-1}c_1(TM)\cdot f_\ast([\Sigma])
+g_{|\Gamma|}-1\geq 0 & \text{if } g_{|\Gamma|}\geq 2.\\
\end{array}
$$
Furthermore, note that in cases (i) and (iii), the complex structure $j$
has an automorphism group of complex dimension $3-k$ and $1$ respectively,
so that in each of these two cases we have a sharper inequality
$$
\begin{array}{ll}
{\em (i^\prime)} \hspace{2mm}|H|^{-1}c_1(TM)\cdot f_\ast([\Sigma])
-1\geq 0 & \text{if } g_{|\Gamma|}=0, k\leq 3,\\
{\em (iii^\prime)}  \hspace{2mm}|H|^{-1}c_1(TM)\cdot f_\ast([\Sigma])
\geq 0 & \text{if } g_{|\Gamma|}=1, k=0.\\
\end{array}
$$

With these inequalities in hand, it now comes to the following observation.

\vspace{2mm}

\noindent{\bf Claim 2}\hspace{2mm}
{\em
$c_1(TM)\cdot C_i\leq 0$ for each of the $J$-holomorphic curves in $\{C_i\}$.
}

\vspace{2mm}

To see this, note that $c_1(K)\cdot C_i=\sum_s n_s C_s\cdot C_i\geq n_i C_i^2$,
so that if $c_1(TM)\cdot C_i>0$, one has $C_i^2\leq c_1(K)\cdot C_i<0$,
and from the adjunction formula, $C_i$ must be an embedded $(-1)$-sphere,
contradicting the minimality assumption on $(M,\omega)$. Hence Claim 2.

Now back to the proof of Lemma 2.4. Fix any $J\in\J_0$, we
consider the $J$-holomorphic curves $\{C_i\}$. For any $C_i$, let $H_i\subset
G$ be the subgroup which leaves $C_i$ invariant, and let 
$f_i:\Sigma_i\rightarrow
M$ be an equivariant $J$-holomorphic map parametrizing $C_i$. Then by Claim 2
there are only two possibilities for $C_i$: (1) $\Gamma_i\equiv\Sigma_i/H_i$
is of genus one and $H_i$ acts on $\Sigma_i$ freely, and $c_1(K)\cdot C_i=0$,
(2) the underlying surface of $\Gamma_i\equiv\Sigma_i/H_i$ has genus
$g_{|\Gamma_i|}\geq 2$, and $c_1(K)\cdot C_i\leq |H_i|\cdot(g_{|\Gamma_i|}-1)$.
Moreover, it follows easily that in case (1), $C_i$ is an embedded
torus with $C_i^2=0$, which is disjoint from the rest of the
$J$-holomorphic curves and is in the complement of the exceptional
orbits.

As for case (2), note that $g_{\Sigma_i}-1\geq |H_i|\cdot (g_{|\Gamma_i|}-1)$
with equality iff $H_i$ acts on $\Sigma_i$ freely. On the other
hand, by the adjunction formula,
$$
2(g_{\Sigma_i}-1)\leq C_i^2+c_1(K)\cdot C_i\leq
\frac{1}{n_i} c_1(K)\cdot C_i+c_1(K)\cdot C_i\leq
(\frac{1}{n_i}+1)\cdot |H_i|\cdot (g_{|\Gamma_i|}-1),
$$
which implies that $n_i=1$, $H_i$ acts freely on $\Sigma_i$, and
$c_1(K)\cdot C_i=C_i^2$. It follows easily that in this case, $C_i$
is embedded with genus $C_i^2+1\geq 2$, disjoint from the rest of the
$J$-holomorphic curves, and lies in the complement of the
exceptional orbits. Lemma 2.4 is thus proved.

\hfill $\Box$

Now we give a proof for parts (b) and (c) of the theorem. First of all, 
observe that the geometric genus $p_g(X)$ is nonzero because the linear 
system $|K_X|$ is nonempty, so that we have $b_2^{+}(X)\geq 2$. Let 
$\omega$ be a K\"{a}hler form on $X$ which is
equivariant under $\text{Aut}(X)^o$. We apply Lemma 2.4 to $(X,\omega)$
with $G=\text{Aut}(X)^o$ (note that $X$ is minimal, and
$b_G^{2,+}=b_2^{+}\geq 2$, so that with Lemma 2.1, the assumptions in
Lemma 2.4 are satisfied), and obtain a set of $J$-holomorphic
curves $\{C_i\}$ as in Lemma 2.4, with
$$
c_1^2=(\sum_i n_i C_i)\cdot (\sum_i n_i C_i)=\sum_i n_i^2 C_i^2.
$$
The fact that $c_1^2>0$ implies that there is at least one $C_i$
with $C_i^2>0$, and since $n_i=1$ for any $C_i$ with $C_i^2>0$, we have
$$
c_1^2=\sum_i C_i^2=\sum_i (\text{ genus}(C_i)-1).
$$
Now recall that $\text{Aut}(X)^o$ operates on $H^2(X;\Q)$
trivially, and $\cup_i C_i$ is invariant under $\text{Aut}(X)^o$.
Thus for any $g\in\text{Aut}(X)^o$, $g\cdot C_i$ is disjoint from $C_i$
if $g\cdot C_i\neq C_i$, which can occur only when $C_i^2=0$. This
implies that $\text{Aut}(X)^o$ leaves $C_i$ invariant for any $C_i$
with $C_i^2>0$. Consequently, for any $C_i$ with $C_i^2>0$, $\text{genus}
(C_i)-1$ is divisible by $|\text{Aut}(X)^o|$, which implies that
$c_1^2$ is divisible by $|\text{Aut}(X)^o|$, and moreover, $\text{genus}(C_i)$
is greater than or equal to the free genus of $\text{Aut}(X)^o$,
which implies part (c) of the theorem.

\vspace{2mm}

{\Small Address: Mathematics Department, Tulane University,
New Orleans, LA 70118, USA.\\
{\it e-mail:} wchen@@math.tulane.edu}

\end{document}